\documentclass[12pt]{article}
\usepackage{amssymb, epsf,chicago, here, amsfonts,latexsym,sym1,graphicx}

\newcommand{\dfrac}[2]{{\textstyle \frac{#1}{#2}}}

\newcommand{\mycite}[1]{{\small \sc \citeNP{#1}}}

\newtheorem{theorem}{Theorem}[section]

\newtheorem{lemma}{Lemma}[section]

\newtheorem{prop}{Proposition}[section]

\newtheorem{defn}{Definition}[section]

\def\argmin{\mbox{argmin}}
\def\a{\alpha}

\def\R{I\!\!R}

\def\d{\delta}

\def\l{\lambda}
\def\labda1{\lambda_1}
\def\labda2{\lambda_2}
\def\e{\epsilon}

\def\s{\sigma}

\def\argmin{\mbox{argmin}}

\def\comment#1{\relax}

\def\=in{\mathop{\rm =}}

\advance \topmargin by -\headheight
\advance \topmargin by -\headsep
\advance \topmargin .2in
\begin{document}
\title{The support reduction algorithm\\ for computing\\ nonparametric
function estimates in mixture models}
\author{Piet Groeneboom, Geurt Jongbloed\thanks{
Research supported by a grant from the Haak Bastiaanse Kuneman
foundation of the Vrije Universiteit.}, Jon A.\ Wellner}
\date{\today}
\maketitle

\begin{abstract}
Vertex direction algorithms have been around for a few decades in
the experimental design and mixture models literature. We briefly
review this type of algorithms and describe a new member of the
family: the support reduction algorithm. The support reduction
algorithm is applied
to the problem of computing nonparametric estimates in two inverse
problems: convex density estimation and
the Gaussian deconvolution problem. Usually,
VD algorithms solve a finite dimensional (version of the) optimization
problem of interest. We introduce a method to solve the true infinite
dimensional optimization problem.
\end{abstract}

\section{Introduction}
\label{sec:intro} \setcounter{equation}{0} During the past decades
emphasis in statistics has shifted from the study of parametric
models to that of semi- or nonparametric models. A big advantage
of these latter models is their flexibility and ability to `let
the data speak for itself'.
However, also problems that were not
usually crucial in the parametric case, turn out to be difficult
in the semiparametric situation. The asymptotic distribution
theory of estimators is one of these problems. The multivariate
central limit theorem and the delta method give the answer to many
questions regarding asymptotic distribution theory in the
parametric setting. For the semiparametric situation, such `basic
tools' are not available. Another problem that is usually easier
to solve in parametric models is the problem of computing
$M$-estimators that are defined as minimizer of a random criterion
function. In a parametric model often estimates can be computed
explicitly or computed using some numerical technique for solving
(low dimensional) convex unconstrained optimization problems like
steepest descent or Newton. In semiparametric models, the
computational issues often boil down to high dimensional
constrained optimization problems.

Apart from algorithms that are known from the general theory of
optimization, algorithms have been designed within the field of
statistics that are particularly useful in certain statistical
applications. Perhaps the best known example of this type is the
{\it Expectation Maximization} (EM) algorithm of \mycite{dlr:77},
that is designed to compute maximum likelihood estimates based on
incomplete data. Another example is the {\it iterative convex
minorant} algorithm that is introduced in \mycite{growe:92} and
further studied in \mycite{jo:98}. That algorithm is based on
techniques known from the theory of isotonic regression as can be
found in \mycite{robwridyk:88} and can be used to compute
nonparametric estimators of distribution functions in
semiparametric models. Another class of algorithms that falls
within this framework is the class of {\it vertex direction} (VD)
algorithms.

In section \ref{sec:vdm} we introduce the general structure
of VD algorithms and mixture models where VD algorithms can
be used to compute nonparametric function estimates. Two
specific examples
of these mixture models will be considered in subsequent
sections: estimating a convex decreasing density and estimating
a mixture of unit variance normal distributions.

In section \ref{sec:simar}, we introduce the support reduction
algorithm as a specific member of the VD family of algorithms.
This algorithm essentially replaces the original infinite
dimensional constrained optimization problem by a sequence of
finite dimensional unconstrained optimization problems. The
algorithm is designed to keep the dimension of these sub-problems
as low as possible. For a specific type of statistical models, the
algorithm seems to be a good candidate to compute sensible
estimators. These are problems that are difficult from the
asymptotic statistical point of view in the sense that the
convergence rate of the estimator is relatively low.

All VD algorithms have to deal with a problem of
minimizing a ``directional derivative'' function over some set of parameters.
There are some variants of these functions. For quadratic objective functions,
we will describe an alternative directional
derivative function in section \ref{sec:quad} that takes more local
information of the objective function into account.

The directional derivative function (and our alternative) are usually nonconvex
functions on a continuum of parameters.
Usually the associated nonconvex minimization problem is circumvented
by considering a fine grid
within the parameter space and minimizing the function only over that
grid. In section \ref{sec:gridless} we propose a method of ``leaving
the grid'', tackling the infinite dimensional optimization problem
rather than the finite dimensional approximation.

Section \ref{sec:lsconv} is devoted to least squares estimation
within a mixture model. The general procedure is given and for the
problem of estimating a convex and decreasing density based on a
sample from it, will be considered in detail. In that situation the
support reduction algorithm boils down to what is called the
iterative cubic spline algorithm in \mycite{grojowe2:01}.

In section \ref{sec:mldecon} the general problem of computing a
maximum likelihood estimate within a mixture model will be
addressed. A Newton procedure based on the support reduction
algorithm will be described. The normal deconvolution problem
will serve as example to illustrate the general approach.

\section{Vertex direction-type algorithms}
\label{sec:vdm} Consider the following type of optimization
problem
\begin{equation}
\label{eq:optprob} \mbox{ minimize } \phi(f) \ \ \mbox{ for }
\ \ f\in C
\end{equation}
where $\phi$ is a convex function defined on (a superset of) a
convex set of functions $C$. We assume throughout that $\phi$
has a unique minimizer over $C$.

\vspace{0.2cm}

\noindent {\bf Assumption A1:} $\phi$ is a convex function on $C$
such that for each $f,g\in
C$ where $\phi$ is finite, the function $t\mapsto \phi(f+t(g-f))$
is continuously differentiable for $t\in (0,1)$.

\vspace{0.2cm} Now define, for each $f\in C$ and $h$ a function
such that for some $\e>0$, $f+\e h\in C$,
$$
D_{\phi}(h;f)=\lim_{\e\downarrow0} \e^{-1}\left(\phi(f+\e
h)-\phi(f)\right)
$$
Note that this quantity exists (possibly equal to $\infty$) by convexity of
$\phi$. As we will see, a choice often made for $h$ is $h=g-f$ for
some arbitrary $g\in C$. In that case we have
$$
D_{\phi}(g-f;f)=\lim_{\e\downarrow0}
\e^{-1}\left(\phi(f+\e(g-f))-\phi(f)\right)
$$

The following simple but important result gives necessary and
sufficient conditions for $\hat{f}$ to be the solution of
(\ref{eq:optprob}).
\begin{lemma}
\label{lem:optcond1} Suppose that $\phi$ satisfies {\bf A1}. Then
$$
\hat{f}=\argmin_{f\in C}\phi(f) \ \ \ \mbox{if and only if} \ \ \
D_{\phi}(g-\hat{f};\hat{f})\ge0 \ \ \ \mbox{for all} \ \  g\in C \, .
$$
\end{lemma}
{\bf Proof:} First we prove $\Rightarrow$. Suppose
$\hat{f}=\argmin_{f\in C}\phi(f)$ and choose $g\in C$ arbitrarily.
Then, for $\e\downarrow 0$
$$
0\le \e^{-1}(\phi(\hat{f}+\e(g-\hat{f}))-\phi(\hat{f}))\downarrow
D_{\phi}(g-\hat{f};\hat{f})
$$
Now $\Leftarrow$. For arbitrary $g\in C$, write $\tau$ for the
convex function $\e\mapsto \phi(\hat{f}+\e(g-\hat{f}))$ and note
that
$$
\phi(g)-\phi(\hat{f})=\tau(1)-\tau(0)\ge
\tau^{\prime}(0+)=D_{\phi}(g-\hat{f};\hat{f})\ge0 \, .
$$
\hfill$\Box$

\vspace{0.2cm}

Consider now the case where $C$ is the convex hull of a class of
functions
\begin{equation}
\label{eq:defF} {\cal F}=\{f_{\theta}\,:\,\theta\in\Theta\subset
\R^k\},
\end{equation}
in the sense that
\begin{equation}
\label{eq:genrepr} C={\rm conv}({\cal
F})=\left\{g=\int_{\Theta}f_{\theta}\,d\mu(\theta)\,: \,\mu \mbox{
probability measure on }\Theta\right\}.
\end{equation}
Here are two examples of {\it mixture models} that fall within this
framework. These examples will reappear in subsequent sections.

\vspace{0.3cm}

\noindent
{\bf Example 1.} {\it (convex decreasing density)}\\
The class of convex decreasing densities on $[0,\infty)$ has
representation (\ref{eq:genrepr}) with
$$
f_{\theta}(x)=\frac{2(\theta-x)}{\theta^2}1_{(0,\theta)}(x),
\,\,\,\, \theta>0 \, .
$$
It is obvious that any (positive) mixture of these
functions is convex and decreasing. Since the mixing measure is
a probability measure, it also follows that the mixture is a
probability density. To see that any convex and decreasing density
can be written as mixture of $f_{\theta}$'s, note that the measure
defined by $d\mu(\theta)=\frac12 \theta^2 df^{\prime}(\theta)$
gives
$$
\int_0^{\infty}f_{\theta}(x)\,d\mu(\theta)=\int_x^{\infty}\frac{2(\theta-x)}{\theta^2}\,d\mu(\theta)=
\int_x^{\infty}(\theta-x)\,df^{\prime}(\theta)=f(x) \, .
$$
Situations where the problem of estimating a convex and decreasing
density based on a sample from it is encountered, can e.g.\ be
found in \mycite{hamp:87} and \mycite{lavsafmei:91}.  \hfill$\Box$

\vspace{0.3cm}

\noindent
{\bf Example 2.} {\it (mixture of unit variance normals)}\\
The {\it Gaussian deconvolution problem} as considered in e.g.\ 
\mycite{growe:92},
entails estimation of a density (and associated mixing distribution)
that belongs to the convex hull of
the class of normal densities with unit variance:
$$
f_{\theta}(x)=\frac1{\sqrt{2\pi}}e^{-\dfrac12(x-\theta)^2} \, .
$$
\hfill$\Box$

\vspace{0.3cm}

   In the examples just considered, usually one has a
sample $X_1,X_2,\ldots,X_n$ from the unknown density
$f \in C$, and wants to
estimate the underlying density $f$ based on that sample. In this
paper we consider two types of nonparametric estimators: least
squares (LS) estimators and maximum likelihood (ML) estimators.

\vspace{0.3cm}

\noindent
{\bf Least Squares estimation}.\\
We define a least squares estimate of the density in $C$ as
minimizer of the function
\begin{equation}
\label{eq:objfunLS} \phi(f)=\frac12\int_0^{\infty}
f(t)^2\,dt-\int_0^{\infty} f(t)\,d\, \FF_n(t)
\end{equation}
over the class $C$. Here $\FF_n$ is the empirical distribution
function of the sample.

The reason for calling this estimator a LS estimator, is the
following heuristic. For any (arbitrary) square integrable density
estimate $\tilde{f}_n$ of $f_0$, one can define the LS estimate as
minimizer of the function
\begin{equation}
\label{eq:LSobj} f\mapsto \frac12\int (f(t)-\tilde{f}_n(t))^2\,dt=
\frac12\int_0^{\infty} f(t)^2\,dt- \int_0^{\infty}
f(t)\tilde{f}_n(t)\,dt+c_{\tilde{f}_n}
\end{equation}
over the class $C$. It is seen that, as far as minimization over
$f$ is concerned, (\ref{eq:LSobj}) only depends on the density
$\tilde{f}_n$ via its distribution function. The objective function
in
(\ref{eq:objfunLS}) is obtained by taking the empirical
distribution function for this estimator, so we take formally
$\tilde{f}_n(t)\,dt=d\, \FF_n(t)$ in (\ref{eq:LSobj}). Note that
for objective function (\ref{eq:LSobj})
$$
D_{\phi}(h;f)=\lim_{\e\downarrow0} \e^{-1}\left(\phi(f+\e
h)-\phi(f)\right)=\int h(x)f(x)\,dx-\int h(x)\, d\, \FF_n(x) \, .
$$

\vspace{0.3cm}

\noindent
{\bf Maximum Likelihood estimation.}\\
As maximum likelihood estimate we define the minimizer of the
function
$$
\phi(f)=-\int \log f(x) \,d\, \FF_n(x)
$$
over the class of densities $C$. Note that for this function
$$
D_{\phi}(h;f)=\lim_{\e\downarrow0} \e^{-1}\left(\phi(f+\e
h)-\phi(f)\right)=-\int\frac{h(x)}{f(x)}\,d\FF_n(x).
$$

\vspace{0.3cm}

\noindent
For both objective functions $\phi$, the function
$D_{\phi}$ has the linearity property stated below.

\vspace{0.3cm}

\noindent {\bf Assumption A2:} the function $\phi$ has the
property that for each $f\in C$ and
$g=\int_{\Theta}f_{\theta}\,d\mu_g(\theta)\in C$,
\begin{equation}
\label{eq:linder}
D_{\phi}(g-f;f)=\int_{\Theta}D_{\phi}(f_{\theta}-f;f)\,d\mu_g(\theta) \, .
\end{equation}

\vspace{0.2cm}

\noindent Under this additional assumption, the nonnegativity
condition in lemma \ref{lem:optcond1} that has to hold for each
$g\in C$, may be restricted to functions $g\in {\cal F}$.
\begin{lemma}
\label{lem:optcondhull} Suppose that $C=$conv$({\cal F})$ with
${\cal F}$ as in (\ref{eq:defF}) and that $\phi$ satisfies {\bf
A1} and {\bf A2}. Then
$$
\hat{f}=\argmin_{f\in C}\phi(f) \ \ \ \mbox{if and only if} \ \ \
D_{\phi}(f_{\theta}-\hat{f};\hat{f})\ge0 \ \ \mbox{for all} \ \
\theta\in\Theta\, .
$$
\end{lemma}
{\bf Proof:} Follows immediately from lemma \ref{lem:optcond1},
the fact that $f_{\theta}\in C$ and (\ref{eq:linder})\hfill$\Box$

\vspace{0.2cm}

For the situation of Lemma \ref{lem:optcondhull}, there is a
variety of algorithms to solve (\ref{eq:optprob}) that can be
called  `of vertex direction (VD) type'. A common feature of VD
algorithms is that they consist of two basic steps. Given a
current iterate $f$, find a value of $\theta$ such that
$D_{\phi}(f_{\theta}-f;f)$ is negative. (If such a value cannot be
found, the current iterate is optimal!) This means that travelling
from the current iterate in the direction of $f_{\theta}$ would
(initially) decrease the value of the function $\phi$.

Having found such a feasible profitable direction from the current
iterate, the next step is to solve some low-dimensional
optimization problem to get to the next iterate.

The original algorithm, proposed by \mycite{fedorov} and
\mycite{wynn} in the context of computing an optimal design, as
well as the algorithm proposed by \mycite{simar:76} (for computing
the maximum likelihood estimate of the mixing distribution in a
Poisson mixture) that we will come back to later, implement the
first step as follows. Given the current $f$, find $\hat{\theta}$
corresponding to the minimizer of $D_{\phi}(f_{\theta}-f;f)$ over
$\Theta$.

\mycite{fedorov} and \mycite{wynn} then propose to take as new
iterate the function
$$
g=(1-\hat{\e})f + \hat{\e}f_{\hat{\theta}}
$$
where $\hat{\e}$ is given by
$$
\hat{\e}=\argmin_{\e\in[0,1]}\phi((1-\e)f+\e f_{\hat{\theta}})\, .
$$
In words, the next iterate is the optimal convex combination of
the current iterate and the most promising vertex in terms of the
directional derivative. It is clear that usually the next iterate
has one more support point than the current iterate.

The {\it vertex exchange algorithm} as proposed in
\mycite{boehning:86}, not only uses the parameter $\hat{\theta}$
corresponding to the minimizer of $D_{\phi}(f_{\theta}-f;f)$, but
also the maximizer $\check{\theta}$ of $D_{\phi}(f_{\theta}-f;f)$
restricted to the support points of the current iterate to get the
direction. Denote by $\mu_f(\{\check{\theta}\})$ the mass assigned
by the mixing distribution corresponding to $f$ to
$\check{\theta}$. Then the direction given by the algorithm is
$f+\mu_f(\{\check{\theta}\})
(f_{\hat{\theta}}-f_{\check{\theta}})$. The new iterate becomes
$$
f+\hat{\e}\mu_f(\{\check{\theta}\})(f_{\hat{\theta}}-f_{\check{\theta}})
$$
where
$$
\hat{\e}=\argmin_{\e\in[0,1]}\phi(f+\e
\mu_f(\{\check{\theta}\})(f_{\hat{\theta}}-f_{\check{\theta}})) \, .
$$
If $\hat{\e}=1$, the point $\check{\theta}$ is eliminated from the
support of the current iterate, and the mass assigned to
$\check{\theta}$ by the `old' mixing distribution, is moved to the
new point $\hat{\theta}$. It is clear that in this algorithm the
number of support points of the iterate can increase by one,
remain the same, but also decrease by one during one iteration
(if $\hat{\e}=1$ and $\check{\theta}$ already belongs to the support).
In specific examples, the number of support points of
the solution $\hat f$ is known to be smaller than a constant $N$ which
only depends on the data (and is known in advance). In the context
of random coefficient regression models, \mycite{mallet:86} proposes
to restrict all iterates to having at most $N$ support points.

Another variation on the theme is due to \mycite{lk:92}. It is
called the {\it intra simplex direction method}. The set of all
local minima $\{\theta_1,\ldots,\theta_m\}$ of
$D_{\phi}(f_{\theta}-f;f)$, where $D_{\phi}$ is negative, is
determined and the optimal convex combination of the current
iterate and all vertices $f_{\theta_1},\ldots,f_{\theta_m}$ is the
new iterate. This final step is to minimize a convex function in
the variables $\e_1.\ldots,\e_m$ under the constraint
$0\le\sum_{i=1}^m\e_i\le1$.

The aforementioned algorithm proposed by \mycite{simar:76} and
further studied in \mycite{boehning:82}, sticks to the original
idea of picking one $\theta$ corresponding to a profitable
direction. The second step differs from those indicated above.
Denote by $S_f$ the set of support points of the mixing measure
corresponding to a function $f\in C$. Then, given $\hat{\theta}$,
the  next iterate is given by
$$
g=\argmin_{h\in C(f)}\phi(h),\ \
\mbox{ where } \ \ C(f)=\{h\in C\,:\,
S_h\subset S_f\cup \{\hat{\theta}\}\} \, .
$$
It is to be noted that support points can (and usually do) vanish
during this second step. Under certain conditions,
\mycite{boehning:82} proves convergence of this algorithm.

In section \ref{sec:simar} we revisit Simar's algorithm and
propose an extension of it that can deal with the case where $C$
is the {\it convex cone} rather than {\it convex hull} generated
by ${\cal F}$. This is convenient for the examples we consider.
Moreover, we will introduce an algorithm that is closely related
to Simar's algorithm: the {\it support reduction algorithm}.

\section{Support reduction and Simar's algorithm}
\label{sec:simar} In Simar's original algorithm,
two optimization problems have to be solved. The first is to minimize
the (usually
nonconvex) function $D_{\phi}(f_{\theta}-f;f)$ in $\theta$. The
second is to minimize $\phi$ over the convex set of functions that
is generated by finitely many functions from ${\cal F}$. In many
examples (including the examples considered here), this  second
step gets more tractable if we were allowed
to minimize over the {\it convex cone} generated by these finitely
many functions in ${\cal F}$. In this section we therefore
consider our function class ${\cal F}$ and the convex cone $C$
generated by it:
$$
C={\rm cone}({\cal
F})=\left\{g=\int_{\Theta}f_{\theta}\,d\mu(\theta)\,: \,\mu \mbox{
positive finite measure on }\Theta\right\} \, .
$$

As will be seen in section \ref{sec:lsconv} and \ref{sec:mldecon},
our two examples fit within this framework of minimizing $\phi$
over the convex cone generated by a set of functions. Assumption
{\bf A2} is now replaced by the following.

\vspace{0.3cm}

\noindent {\bf Assumption A2$^\prime$:} the function $\phi$ has the
property that for each $f\in C={\rm cone}({\cal F})$ and
$g=\int_{\Theta}f_{\theta}\,d\mu_g(\theta)\in C$,
\begin{equation}
\label{eq:linder2} D_{\phi}(g;f)=\int
D_{\phi}(f_{\theta};f)\,d\mu_g(\theta) \, .
\end{equation}

\vspace{0.2cm}

\noindent {\bf Remark}. Suppose that $h_1$ and $h_2$ are such that
for a small positive $\e$, $f+\e h_i\in C$ for $i=1,2$. Then,
since $C$ is convex, we have that  $f+\dfrac12\e (h_1+h_2)\in C$,
and $D_{\phi}(\cdot;f)$ is well defined at $h_1$, $h_2$ and
$h_1+h_2$. Assumption {\bf A2$^\prime$} then implies the following
linearity property:
\begin{equation}
\label{eq:linpropder} D_{\phi}(h_1+h_2;f)=\int
D_{\phi}(f_{\theta};f)\,d\mu_{h_1+h_2}(\theta)=
\int D_{\phi}(f_{\theta};f)\,d\left(\mu_{h_1}+\mu_{h_2}\right)(\theta)
=D_{\phi}(h_1;f)+D_{\phi}(h_2;f) \, .
\end{equation}

\vspace{0.2cm}

\noindent {\bf Remark}. Assumption {\bf A2$^\prime$} implies {\bf A2} for
$g\in{\rm conv}({\cal F})$. Indeed, take
$g=\int_{\Theta}f_{\theta}\,d\mu_g(\theta)\in {\rm conv}({\cal
F})$, meaning that $\mu_g$ is a probability measure. Then we have,
also using (\ref{eq:linpropder}),
\begin{eqnarray*}
D_{\phi}(g-f;f)&=&D_{\phi}(g;f)-D_{\phi}(f;f)
=\int D_{\phi}(f_{\theta};f)\,d\mu_g(\theta)-D_{\phi}(f;f) \\
&=&\int D_{\phi}(f_{\theta};f)-D_{\phi}(f;f)\,d\mu_g(\theta) =\int
D_{\phi}(f_{\theta}-f;f)\,d\mu_g(\theta) \, .
\end{eqnarray*}

Let us formulate a result for a generated cone analogous to lemma
\ref{lem:optcondhull}.
\begin{lemma}
\label{lem:optcondcone} Let $C={\rm cone}({\cal F})$ and $\phi$
satisfy {\bf A1} and {\bf A2$^\prime$}. Suppose that the measure
$\mu_{\hat{f}}$ in
$\hat{f}=\int_{\Theta}f_{\theta}\,d\mu_{\hat{f}}(\theta)$ has
finite support. Then
\begin{equation}
\label{eq:lemstat} \hat{f}=\argmin_{f\in C}\phi(f) \ \ \
\mbox{if and only if} \ \ \
D_{\phi}(f_{\theta};\hat{f})\left\{
\begin{array}{ll}
\ge0 & \mbox{for all } \  \theta\in\Theta\\
=0   & \mbox{for all } \  \theta\in\mbox{supp}(\mu_{\hat{f}}) \, .
\end{array}
\right.
\end{equation}
\end{lemma}
{\bf Proof:} If $\hat{f}=\argmin_{f\in C}\phi(f)$, then we have by
{\bf A1} that
$$
D_{\phi}(\hat{f};\hat{f})=\lim_{\e\to0}\e^{-1}
\left(\phi((1+\e)\hat{f})-\phi(\hat{f})\right)=0 \, .
$$
Hence, by (\ref{eq:linpropder}) and lemma \ref{lem:optcond1}, we
have for all $\theta\in\Theta$
\begin{equation}
\label{eq:decomp}
D_{\phi}(f_{\theta};\hat{f})=D_{\phi}(f_{\theta}-\hat{f}+\hat{f};\hat{f})=
D_{\phi}(f_{\theta}-\hat{f};\hat{f})+D_{\phi}(\hat{f};\hat{f})=
D_{\phi}(f_{\theta}-\hat{f};\hat{f})\ge 0
\end{equation}
In view of property (\ref{eq:linder2}), we have
$$
0=D_{\phi}(\hat{f};\hat{f})=\int D_{\phi}(f_{\theta};\hat{f}) \,
d\mu_{\hat{f}}(\theta) \, .
$$
In the presence of the inequalities in (\ref{eq:decomp}) we
therefore have that $D_{\phi}(f_{\theta};\hat{f})=0$ on the
support of $\mu_{\hat{f}}$ necessarily.

Conversely, if $\hat{f}$ satisfies the (in)equalities given in
(\ref{eq:lemstat}) above, we have for any $f\in C$ that
$$
\phi(f)-\phi(\hat{f})\ge
D_{\phi}(f-\hat{f};\hat{f})=D_{\phi}(f;\hat{f})=\int
D_{\phi}(f_{\theta};\hat{f}) \, d\mu_{f}(\theta)\ge0 \, .
$$
\hfill$\Box$

\vspace{0.2cm}

\noindent {\bf Remark.} The assumption that the support of
$\mu_{\hat{f}}$ is finite seems to be restrictive and unnatural.
However, there are many problems (including our examples) where
this is true. Of course, if $\Theta$ is finite it is trivially
true (this e.g.\ covers interval censoring problems). Moreover,
maximum likelihood estimators in mixture models usually have this
property (\mycite{lindsay:95}, theorem 18, section 5.2).

\vspace{0.2cm}

\noindent Below we give the pseudo code for Simar's algorithm
constructed for a cone and also for the support reduction
algorithm we propose. In fact, as will be seen below, the support
reduction algorithm is Simar's algorithm where one substep is not
completely followed till the end.

\noindent \hspace{0.8cm}

\hrule

\hspace{0.3cm}

\large \noindent {\bf Basic Simar- and support reduction algorithm
for a cone}

\normalsize

\hspace{0.3cm}

\hrule

\hspace{0.3cm}

\noindent
{\bf Input:}\\
$\eta>0$: accuracy parameter;\\
$\theta^{(0)}\in\Theta$: starting value;\\
$f=\argmin_{g\in C\,:\,S_g=\{\theta^{(0)}\}}\phi(g)$;

\begin{tabbing}
{\bf be}\={\bf gin}\\
          \>{\bf while}
$\min_{\theta\in\Theta}D_{\phi}(f_{\theta};f)<-\eta$
              {\bf do}\\
          \>{\bf be}\={\bf gin}\\
          \>        \>$\hat{\theta}:=\argmin_{\theta\in\Theta}
                        D_{\phi}(f_{\theta};f)$;\\
          \>        \>$S^*:=S_f\cup\{\hat{\theta}\}$;\\
          \>        \>$f:=\argmin_{g\in C\,:\,S_g\subset S^*}\phi(g)$;
({\it Simar})\\
          \>        \>$f:=\argmin_{g\in C\,:\,S_g\subset_s S^*}\phi(g)$;
({\it Support
reduction})\\
          \>{\bf end};\\
{\bf end.}
\end{tabbing}
\hrule

\vspace{0.5cm}

\noindent The meaning of `$\subset_s$' will become clear in the
sequel. For both algorithms, there are two finite dimensional
optimizations that have to be performed. The first one is over
$\Theta$. In general the function $\theta\mapsto
D_{\phi}(f_{\theta};f)$ is nonconvex and minimizing such a
function is usually difficult. Hence, in each setting one should
try to take advantage of the specific features of that problem to
attack this first optimization problem. Usually one can restrict
the minimization to a bounded subset of $\Theta$ and use
a fine (finite) grid in this subset instead of the whole set $\Theta$.
Then the minimization reduces to finding the minimal element in a
(long) vector. After that, it is possible to `leave the grid' in a way
as described in section \ref{sec:gridless}. Sometimes (e.g.\ when
computing the ML estimator of a distribution function based on
interval censored observations) it is even possible to select a
finite subset of $\Theta$, based on the data, such that the minimizer of
$\phi$ over $C$ is contained in the convex hull of the corresponding
finitely many generators.
In subsequent sections, we will address this matter more specifically in the
examples.

The second optimization in the algorithm is over a convex cone
that is spanned by finitely many functions $f_{\theta}$ in ${\cal
F}$. Lemma \ref{lem:optcondcone} gives a characterization of such
a function (applying the lemma to the finite subset $S^*$ of
$\Theta$ instead of $\Theta$ itself). We propose the following
general way of solving this finite dimensional constrained
optimization problem in Simar's algorithm. In passing it will
become clear what the support reduction algorithm does.

Given the current iterate and the new support point
$\hat{\theta}$, consider the {\it linear space} $L$ spanned by the
finitely many functions $\{f_{\theta}\,:\,\theta\in\ S^*\}$:
$$
L_{S^*}=\left\{ g=\int_{\Theta} f_{\theta}\,d\s(\theta)\,:\,\s \mbox{ is
a finite signed measure on }S^*\right\} \, ,
$$
and determine
$$
f_{u}^{(0)}=\argmin_{g\in L_{S^*}}\phi(g)=\int
f_{\theta}\,d\s_{f_{u,S^*}}(\theta)= \sum_{\theta\in
S^*}f_{\theta}\s_{f_{u,S^*}}(\{\theta\}) \, .
$$
We assume $\phi$ has a smooth convex extension to the space $L_{S^*}$.
In our examples and many others this is certainly the case. This
optimization corresponds to finding a solution of a finite system of
equations. Of course, $f_{u}^{(0)}$ will in general not be an
element of $C$, since certain coefficients
$\s_{f_{u}^{(0)}}(\{\theta\})$ may be negative. Nevertheless we
can {\it always} move from $f$ towards $f_{u}^{(0)}$ and stay
within the class $C$ initially. This is a consequence of the fact
that the coefficient $\s_{f_{u}^{(0)}}(\{\hat{\theta}\})$ of
$f_{\hat{\theta}}$ in $f_{u}^{(0)}$ will be {\it strictly}
positive. Indeed,
$$
0>\lim_{\e\downarrow0}\e^{-1}\left(\phi(f+\e
f_{u}^{(0)})-\phi(f)\right) =\int
D_{\phi}(f_{\theta};f)\,d\sigma_{f_{u}^{(0)}}(\theta)=
\sigma_{f_{u}^{(0)}}(\{\hat{\theta}\})D_{\phi}(f_{\hat{\theta}};f)
$$
and $D_{\phi}(f_{\hat{\theta}};f)<0$ by choice of $\hat{\theta}$.
If $f_u^{(0)}\in C$ then take this as next iterate. Otherwise
define
\begin{equation}
\label{eq:maxmin}
\hat{\lambda}=\max\{\l\in(0,1]\,:\,f+\l(f_{u}^{(0)}-f)\in C\}=
\min_{\theta\in S^*\,:\,\s_{f_{u,S^*}}(\{\theta\})<0}
(1-\s_{f_{u}^{(0)}}(\{\theta\})\slash\s_{f}(\{\theta\}))^{-1}
\end{equation}
and take as next iterate the function $f+\hat{\l}(f_{u}^{(0)}-f)$
and delete the support point $\check{\theta}\in S^*$ where the
minimum in the expression on the right hand side of
(\ref{eq:maxmin}) is attained from the support set:
$$
S^{*(1)}=S^{*}\setminus \{\check{\theta}\} \, .
$$
Then compute the next unrestricted minimizer
$$
f_{u}^{(1)}=\argmin_{g\in L_{S^{*(1)}}}\phi(g) \, .
$$
If this function differs from the current iterate, again a step of
positive length can be made in this direction, since for all
$\theta\in S^{*(1)}$, $\s_{f_{u,S^*}^{(1)}}(\{\theta\})>0$. If we
can go all the way to $f_{u,S^*}^{(1)}$, stop the iteration, and
else delete the support point as it was done in the first step.
This deletion of support points can be continued until we get a
subset $S^{*(j)}\subset S^*$ and a function  $f_{u,S^*}^{(j)}\in
C$ with support set $S^{*(j)}$ such that
$$
D_{\phi}(f_{\theta};f_{u,S^*}^{(j)})=0 \,\,\mbox{for all } \
\theta\in S^{*(j)} \, .
$$
The specific set $S^{*(j)}$ obtained in this way as subset of
$S^*$ is denoted by `$\subset_s$', and this gives the next iterate
in the support reduction algorithm. Note that the function $\phi$
is decreased all the way during the iterations of this substep.

For Simar's algorithm, one should check for the points in
$S^*\setminus S^{*(j)}$ whether the value $\phi$ can be improved
upon by adding such points to the current support.
The natural thing to do then is to take the value of
$\theta\in S^{*(j)}$ where $D_{\phi}(f_{\theta};f)$ is minimal and
add this to the support. In the support reduction algorithm we
skip the adding of deleted points from $S^*$ and allow the next
support point to be chosen without restriction from the whole set $\Theta$.

Let us summarize the steps sketched above to determine
$f:=\argmin_{g\in C\,:\,S_g\subset_s S^*}\phi(g)$ in pseudo code.

\noindent \hspace{0.8cm}

\hrule

\hspace{0.3cm}

\large \noindent {\bf Support reduction step}

\normalsize

\hspace{0.3cm}

\hrule

\hspace{0.3cm}

\noindent
{\bf Input:}\\
$f^{(0)}=f\in C$: minimizer of $\phi$ over subset of $C$
consisting of functions with same
support $S_f$;\\
$S^{*(0)}=S^*=S_f\cup \{\hat{\theta}\}$: finite set of support points;\\
$j:=0$;\\
\begin{tabbing}
{\bf be}\={\bf gin}\\
          \> $f_{u}^{(j)}=\argmin_{g\in L_{S^{*(j)}}}\phi(g)$;\\
          \>{\bf while} $f_{u}^{(j)}\not\in C$ {\bf do}\\
          \>{\bf be}\={\bf gin}\\
          \>        \>$j:=j+1$;\\
          \>        \>$\check{\Theta}=\{\theta\in
S^*\,:\,\s_{f_{u}^{(j-1)}}(\{\theta\})<0$
and
$\s_{f_{u}^{(j-1)}}(\{\theta\})\slash\s_{f^{(j-1)}}(\{\theta\})$
is minimal$\}$;\\
          \>
\>$\hat{\lambda}=(1-\s_{f_{u}^{(j-1)}}(\{\theta\})\slash\s_{f^{(j-1)}}(\{\theta\}))^{-1}$ 

for some $\theta\in\check{\Theta}$;\\

          \>
\>$f^{(j)}=f^{(j-1)}+\hat{\lambda}(f_{u}^{(j-1)}-f^{(j-1)})$;\\
          \>        \>$S^{*(j)}:=S^{*(j-1)}\setminus\check{\Theta}$;\\

          \>        \>$f_{u}^{(j)}=\argmin_{g\in
L_{S^{*(j)}}}\phi(g)$;\\
          \>{\bf end};\\
          \>$f:=f_{u}^{(j)}\in C$: minimizer of $\phi$ over subset of $C$
                   consisting of\\
          \>  \phantom{$f:=f_{u}^{(j)}\in C$: min} functions with same
support
$S_f=S^{*(j)}\subset_s S^*$;\\
{\bf end.}
\end{tabbing}
\hrule

\vspace{0.5cm}

\noindent We now see that the basic building stone of the
algorithm is an {\it unrestricted minimization} of the function
$\phi$. As we will see in the sections \ref{sec:lsconv} and
\ref{sec:mldecon}, there are efficient algorithms to solve this
kind of optimization problems in specific situations.

Before applying the algorithm to concrete problems, let us
consider the convergence issue. The theorem below (the proof of
which is inspired by \mycite{boehning:82}) states that the
algorithms considered in this section indeed converge to the
solution of the optimization problem. To get this, we need one
additional condition on the function $\phi$. This condition is
needed to guarantee that a strictly negative value of
$D_{\phi}(f_{\theta};f)$ for some $\theta$ means that the next
iterate will have some minimal decrease in $\phi$-value.

\vspace{0.2cm}

\noindent {\bf Assumption A3:} For any specific starting function
$f^{(0)}\in C$  with $\phi(f^{(0)})<\infty$, there exists an
$\bar{\e}\in(0,1]$ such that for all $f\in C$ with
$\phi(f)<\phi(f^{(0)})$ and $\theta\in\Theta$, the following
implication holds:
$$
D_{\phi}(f_{\theta}-f;f)\le -\delta<0\Rightarrow
\phi(f+\e(f_{\theta}-f))-\phi(f)\le -\dfrac12\e\d \,\,\ \
\mbox{for all} \ \
\e\in (0,\bar{\e}]
$$

\vspace{0.2cm}

\noindent We will see that this assumption holds for the problems
we will address in subsequent sections.

\vspace{0.2cm}
\begin{theorem}
\label{th:convproof} Denote by $f_n$ a sequence generated by one
of the algorithms introduced here. Then, under the assumptions
{\bf A1}, {\bf A2$^\prime$} and {\bf A3} we have that $\phi(f_n)\downarrow
\phi(\hat{f})$ as $n\to\infty$.
\end{theorem}
{\bf Proof:} Since we have for each $n$ that
$$
f_n=\argmin_{f\in C\,:\ S_f=S_{f_n}}\phi(f),
$$
we have by assumption {\bf A1} that $D_{\phi}(f_n;f_n)=0$. Hence,
by (\ref{eq:linpropder}), we have for all $n\ge0$
$$
D_{\phi}(f_{\theta}-f_n;f_n)=D_{\phi}(f_{\theta};f_n)\,\,
\ \ \mbox{for all} \ \ \theta\in\Theta \, .
$$
Since $\phi(f_n)$ is a bounded and decreasing sequence of real
numbers, it decreases to a limit. Assume for the moment that
$\phi(f_n)\downarrow \phi^*=\phi(\hat{f})+\d>\phi(\hat{f})$ for
some $\d>0$. We will extract a contradiction.

Take $\theta_n$ such that
$D_{\phi}(f_{\theta_n};f_n)\le\frac12\inf_{\theta\in\Theta}
D_{\phi}(f_{\theta};f_n)$. Then we get
\begin{eqnarray}
D_{\phi}(f_{\theta_n}-f_n;f_n)&=&D_{\phi}(f_{\theta_n};f_n)\le\frac12\inf_{\theta\in\Theta}
D_{\phi}(f_{\theta};f_n)\le
\int \frac12D_{\phi}(f_{\theta};f_n)\,d\mu_{\hat{f}}(\theta)\nonumber\\
&=&\frac12D_{\phi}(\hat{f}-f_n;f_n)\le\frac12\left(\phi(\hat{f})-\phi(f_n)\right)\le
\frac12\left(\phi(\hat{f})-\phi^*\right)= -\delta/2 \, .
\label{eq:ineqseq}
\end{eqnarray}
Again by monotonicity of $\phi(f_n)$, we have that $\phi(f_n)\le
\phi(f^{(0)})$ for all $n$, and assumption {\bf A3} gives
\begin{equation}
\label{eq:finineq} \phi(f_{n+1})\le
\phi(f_n+\bar{\e}(f_{\theta_n}-f_n))\le
\phi(f_n)-\dfrac14\bar{\e}\d \ \ \ \mbox{for all} \ \  n\, .
\end{equation}
This contradicts the fact that $\phi(f_n)$ converges. \hfill$\Box$

\vspace{0.2cm}

\noindent In view of the convergence proof, there are some
adaptations of the algorithms that do not destroy the convergence
property of the algorithm. The first adaptation has to do with the
choice of the most promising vertex. If the function $D_{\phi}$ on
$\Theta$ is replaced by a function
$$
\tilde{D}_{\phi}(f_{\theta};f)=w(\theta)D_{\phi}(f_{\theta};f)
$$
where $w$ is some strictly positive weight function on $\Theta$
such that
$$
0<\underline{w}\le w(\theta)\le \overline{w}<\infty
\ \ \mbox{ for all }\ \ \theta\in\Theta\,.
$$
Equation (\ref{eq:ineqseq}) would then change to
$$
D_{\phi}(f_{\theta_n}-f_n;f_n)=D_{\phi}(f_{\theta_n};f_n)\le
(\underline{w}\slash \overline{w})\inf_{\theta\in\Theta}
D_{\phi}(f_{\theta};f_n)\le
   -(\underline{w}\slash \overline{w})\delta \, ,
$$
and the argument goes through with $\delta$ replaced by
$(\underline{w}\slash \overline{w})\delta$. Similarly, {\bf A3}
will also hold for $\tilde{D}_{\phi}$ if it holds for $D_{\phi}$.
In section \ref{sec:quad} we will use this idea to define an
alternative directional derivative function.

The second adaptation is the following.  If it is possible after
{\it reduction by deletion} of support points to do an extra step
of {\it reduction by replacing two support points by a third} or
move a support point slightly without increasing the function
$\phi$, this will not prevent the algorithm from converging. This
immediately follows from (\ref{eq:finineq}). Indeed, if we replace
the iterate $f_n$ that would be obtained by the original method by
$\tilde{f}_n$ which satisfies
$$
\tilde{f}_n=\argmin_{f\in C\,:\ S_f=S_{\tilde{f}_n}}\phi(f)\mbox{
and } \phi(\tilde{f}_n)\le\phi(f_n),
$$
the inequality (\ref{eq:finineq}) also holds for $\tilde{f}_{n+1}$
instead of $f_{n+1}$ and the proof goes through. This adaptation of
the algorithm will be discussed more elaborately in section \ref{sec:gridless}.

\section{Alternative directional derivative}
\label{sec:quad}
Consider a quadratic objective function $\phi_q$ on
$C={\rm cone}(\{f_{\theta}\,:\,\theta\in\Theta\})$. The objective function
in the LS estimation context is quadratic automatically and in
section \ref{sec:mldecon}
we will use a Newton algorithm to solve the ML estimation problem. In
that algorithm
a quadratic approximation of the objective function is minimized
during each iteration.

The
function $\phi_q$ is quadratic in $f$. Hence, along line segments
in the linear space spanned by the functions $f_{\theta_1},\ldots
f_{\theta_p}$, the function is also quadratic as a function of one
variable. Along such segments (or lines), the function $\phi_q$
can therefore be minimized explicitly. Given a `current iterate'
$g$ in the algorithm, we consider for each $\theta\in \Theta$ the
following function (alternative choice is to take $f_{\theta}-g$
instead of $f_{\theta}$):
$$
\e\mapsto\phi_q(g+\e
f_{\theta})-\phi_q(g)=c_1(\theta,g)\e+\frac12\e^2c_2(\theta) \, .
$$
Typically, $c_2>0$, so that
$$
\hat{\e}=\hat{\e}_{\theta}=\argmin_{\e}\phi_q(g+\e
f_{\theta})=-\frac{c_1(\theta,g)}{c_2(\theta)}
$$
is the optimal move along the line connecting $g$ and
$g+f_{\theta}$.

In order to have descent direction, we only  consider points
$\theta$ where $c_1(\theta,g)<0$. In that case, $\hat{\e}>0$. As new vertex,
we then define
$$
\hat{\theta}=\argmin_{\theta\in
\Theta\,:\,c_1(\theta)<0}\phi_q(g+\hat{\e}(\theta)f_{\theta})
=\argmin_{\theta\in
\Theta\,:\,c_1(\theta)<0}-\frac{c_1(\theta,g)^2}{2c_2(\theta)}=
\argmin_{\theta\in S} \frac{c_1(\theta,g)}{\sqrt{c_2(\theta)}} \, .
$$

\section{A `gridless' implementation}
\label{sec:gridless}
For a practical implementation of the step of
selecting a new support point, we propose to fix a fine grid
$\Theta_{\delta}$ in $\Theta$ and run the whole algorithm with
$\Theta_{\delta}$ instead of $\Theta$. Having a precise
approximation of the minimizer $f_{grid}$ of $\phi$ over this
finite dimensional cone, one can make the algorithm `gridless' by
fine tuning in the support points. This can be done by augmenting
a step at each iteration in the spirit of the second remark after
theorem \ref{th:convproof}.

Write $f$ for the current iterate (at the first `fine tuning step',
this is $f_{grid}$) and define
$$
\tau(h_1,h_2,\ldots,h_p)=\tau(h_1,h_2,\ldots,h_p;f)
=\phi\left(\sum_{i=1}^p\alpha_i
f_{\theta_i+h_i}\right)-\phi\left(\sum_{i=1}^p\alpha_i
f_{\theta_i}\right)
$$
with $\alpha_1,\ldots,\alpha_p$ fixed and $h=(h_1,\ldots,h_p)^T$
varying over a neighborhood of zero in $\R^p$. The function $\tau$
represents the value of the objective function if the masses
$\alpha_i$ are fixed and
the current support points are shifted a bit. Abusing
notation slightly, write
$$
f_h=\sum_{i=1}^p\alpha_i f_{\theta_i+h_i}
$$
and note
that for the least squares objective function (under mild
smoothness assumptions on the parameterization of $f_{\theta}$)
\begin{equation}
\label{eq:nabtau1}
\frac{\partial\tau}{\partial h_i}(h_1,h_2,\ldots,h_p)=\alpha_i\int
\dot{f}_{\theta_i+h_i}(x)f_h(x)\,dx-\alpha_i
\int\dot{f}_{\theta_i+h_i}(x)\,d\FF_n(x)
\end{equation}
and for the maximum likelihood objective function
\begin{equation}
\label{eq:nabtau2}
\frac{\partial\tau}{\partial h_i}(h_1,h_2,\ldots,h_p)=-\alpha_i\int
\frac{\dot{f}_{\theta_i+h_i}(x)}{f_h(x)} \,d\, \FF_n(x) \, .
\end{equation}
In particular, note that
$$
\frac{\partial\tau}{\partial h_i}(0)=\alpha_i
\frac{d}{d\theta}D_{\phi}(f_{\theta};f)|_{\theta=\theta_i}
$$
for both objective functions.
Hence, the partial derivatives of $\tau$ at zero are visualized in the graph
of $\theta\mapsto D_{\phi}(f_{\theta};f)$ for both objective functions.
Qualitatively, the interpretation of the partial derivatives of
$\tau$ is that if $\frac{\partial}{\partial h_i}\tau(0)<0$,
shifting the support point $\theta_i$ slightly to the right
will result in a decrease of the
objective function. For the moment, fix $h\in \R^p$ with
$\|h\|_2=1$ and consider the function
$$
\mu_h(\epsilon)=\tau(\epsilon h)
$$
on an interval $[0,\epsilon_0]$ for some small $\epsilon_0>0$.
Note that $\tau_h(0)=0$. Then (again under mild smoothness
assumptions) the derivative of $\mu_h$
is given by
$$
\mu_h^{\prime}(\epsilon)=h^T\nabla\tau(\epsilon h)
$$
where $\nabla\tau(\epsilon h)$ is given either by
(\ref{eq:nabtau1})
or (\ref{eq:nabtau2}), depending on the objective function.
Taking $\epsilon=0$, we see that the `most promising' direction to move,
is the direction $-\nabla\tau(0)$, the direction of steepest descent. From
now on take this direction. The aim is now to move the support
points in this direction to get a sufficient decrease in the
objective function. This means that $\mu_h$ is to be minimized as
a function of $\epsilon\in(0,\epsilon_0)$, or at least a value
$\epsilon$ has to be determined such that $\mu_h(\epsilon)$ is
negative. Note that the function $\mu_h$ is {\it nonconvex} in
general. We determine the step length by the method of regula
falsi on the derivative $\mu_h$. At zero this function is zero.
Define $\e_l=0$ and $\e_u=\e_0$. If $\mu_h(\e_u)<0$ then take this
$\e=\e_0$. Otherwise proceed as follows.
$$
\e_n=\frac{\e_l\mu_h^{\prime}(\e_u)-\e_u\mu_h^{\prime}(\e_l)}
{\mu_h^{\prime}(\e_u)-\mu_h^{\prime}(\e_l)} \, .
$$
If $\mu_h^{\prime}(\e_n)>0$ define $\e_u=\e_n$ whereas if
$\mu_h^{\prime}(\e_n)<0$ define $\e_l=\e_n$. This process can be
iterated till $\mu_h^{\prime}(\e_n)$ is sufficiently small in
absolute value. This regula falsi method comes up with a
stationary point of $\mu_h$. If the $\mu_h(\e_n)$ is positive, the
line search procedure should be repeated with $\e_0=c \e_n$ for
some $0<c<1$ (usually close to one). In our experience this step
is hardly ever necessary, but conceptually it is needed. The
procedure will (in case $\e\neq \e_0$) lead to a stationary point
of $\mu_h$ with $\mu_h(\e)<\mu_h(0)=0$. Actually, $\e$ will
usually correspond to the smallest local minimum of $\mu_h$.

Next, define
$$
\bar f:=\sum_{i=1}^p\alpha_if_{\theta_i+\e h_i} \, .
$$
The new iterate $f$ is finally obtained by minimizing $\phi$ over
the cone generated by $\{f_{\theta_i+\e h_i}\,:\,1\le i\le p\}$.
This function satisfies the conditions needed at the beginning of
the just described `fine tuning' step. Hence, it can be iterated
till the norm of $\mu_h^{\prime}(0)$ is sufficiently small.

\section{LS estimation of a convex density}
\label{sec:lsconv} In this section we study the problem of
computing the least squares estimate of a convex and decreasing
density on $[0,\infty)$. In  \mycite{grojowe2:01}, it is shown
that the (uniquely defined) minimizer of the convex function
$\phi$ over conv(${\cal F}$) is the same as the minimizer of
$\phi$ over cone(${\cal F}$). It is also shown that there only
functions $f_{\theta}$ with $\theta\in [x_1,K]$ for some
$K<\infty$ have to be considered in the optimization, since the
optimal function has no change of slope at a location to the left
of $x_1$ and has compact support. Hence we are in the situation of
section \ref{sec:simar}. Moreover, we have
$$
D_{\phi}(f_{\theta};f)=\int_0^{\infty}
f_{\theta}(x)f(x)\,dx-\int_0^{\infty} f_{\theta}(x)\,d\FF_n(x)=
\frac{2}{\theta^2}\left(H(\theta;f)-Y_n(\theta)\right)
$$
where
$$
H(\theta;f)=\int_{x=0}^{\theta}\int_{0}^{x}f(y)\,dy\,dx \mbox{ and
}Y_n(\theta)=\int_0^{\theta}\FF_n(x)\,dx \, ;
$$
here we use the same notation as in \mycite{grojowe2:01}. Note
that the assumptions {\bf A1}, {\bf A2$^{\prime}$} and {\bf A3} are satisfied
in this situation. For {\bf A3} note that
\begin{equation}
\label{eq:exactexpan} \phi(f+\e(f_{\theta}-f))=\phi(f)+\e
D_{\phi}(f_{\theta}-f;f)+ \dfrac12
\e^2\int_0^{\infty}(f_{\theta}(x)-f(x))^2\,dx
\end{equation}
and that for $\theta\in [x_1,K]$
$$
\int_0^{\infty}(f_{\theta}(x)-f(x))^2\,dx=\frac{4}{3\theta}-\frac{4}{\theta^2}H(\theta;f)
+\int_0^{\infty}f(x)^2\,dx\le M
$$
for some big finite $M$ not depending on $\theta$.

Let us now consider the support reduction algorithm.
To start this algorithm,
we fix a starting value $\theta^{(0)}>x_n$.
Then we determine the function
$c f_{\theta_0}$ minimizing $\phi$ as function of $c>0$. To
this end we need the value $c$ that minimizes
$$
c\mapsto
\phi(cf_{\theta^{(0)}})=\dfrac12c^2\int_0^{\infty}f_{\theta^{(0)}}(x)^2\,dx-
c\int_0^{\infty}f_{\theta^{(0)}}(x)\,d\FF_n(x)=
\frac{2c^2}{3\theta^{(0)}}-\frac{2c(\theta^{(0)}-\bar{x}_n)}{(\theta^{(0)})^2}
$$
giving $c=\dfrac32(1-\bar{x}_n/\theta^{(0)})$. If $x_n<3\bar{x}_n$,
one could also
choose to take $\theta^{(0)}=3\bar{x}_n$, so that the starting
function $f^{(0)}$ would
be a density.

The two main
steps are minimizing $D_{\phi}(f_{\theta};f)$ as a function of
$\theta$ and minimizing the function $\phi$ over the space of
piecewise linear functions with bend points in a finite set $S^*$.
For the first step, we follow the line of thought given in section
\ref{sec:quad}.
In this example we have for all $\epsilon>0$ that
$$
\phi(f+\e f_{\theta})=\phi(f)+ \e c_1(\theta,f)+\frac12 \e^2 c_2(\theta)
\mbox{ with }
c_1(\theta,f)=D_{\phi}(f_{\theta};f)
\ \ \mbox{ and } \ \ c_2(\theta)=\frac{4}{3\theta} \, .
$$
Hence, we have as `alternative directional derivative'
function
$$
\tilde{D}_{\phi}(f_{\theta};f)=\frac{c_1(\theta,f)}{\sqrt{c_2(\theta)}}
\simeq \sqrt{\theta}D_{\phi}(f_{\theta};f)
$$
where $\simeq$ denotes `equality apart from a positive multiplicative
constant'.
Note
that, since $w(\theta)=\sqrt{\theta}$ is strictly positive and
uniformly bounded away from zero and infinity on $\Theta=[x_1,K]$,
we are in the situation described below theorem
\ref{th:convproof}. Note that
$\theta\mapsto\tilde{D}_{\phi}(f_{\theta};f)$ is continuous,
$$
\tilde{D}_{\phi}(f_{\theta};f)=0 \ \
\mbox{at} \ \  \theta=0 \ \ \mbox{ and } \ \
\lim_{\theta\to\infty}\tilde{D}_{\phi}(f_{\theta};f)=0 \, .
$$
Hence, if $\tilde{D}_{\phi}(f_{\theta};f)<0$ for some $\theta$, it
attains its minimal value.

The second step in the algorithm boils down to the following
procedure. Write $S^*=S_f\cup
\{\hat{\theta}\}=\{ \theta_1,\theta_2,\ldots,\theta_m\}$ with
$\theta_1<\cdots<\theta_m$ and construct a cubic spline $P$ with
knots $\{\theta_1,\theta_2,\ldots,\theta_m\}$ such that
\begin{equation}
\label{eq:spline} P(\theta)=Y_n(\theta) \ \ \
\mbox{for all} \ \  \theta\in S^*, \ \
P(0)=P'(0)=P''(\theta_m)=0 \, .
\end{equation}
Note that the second derivative of this cubic spline minimizes the
function $\phi$ within the class of linear splines $l$ with knots
concentrated on the set $\{\theta_1,\theta_2,\ldots,\theta_m\}$
subject to the boundary constraint that $l(\theta_m)=0$. This
follows by setting the derivatives of $\phi$ in the directions
$f_{\theta_j}$ equal to zero, i.e.\ $D_{\phi}(f_{\theta_j};f)=0$.

Figure 1 shows the results of the SR algorithm based on a sample
of size $500$ from the standard exponential distribution. The solution
on an equidistant grid in $[0,3x_{(n)}]=[0,16.5]$ consisting of $1000$ points
was obtained after 33 iterations. Furthermore, we used accuracy parameter
$\eta=10^{-10}$.

\vspace{0.2cm}

\begin{figure}[!ht]
\begin{center}
\leavevmode
\includegraphics[width=12cm]{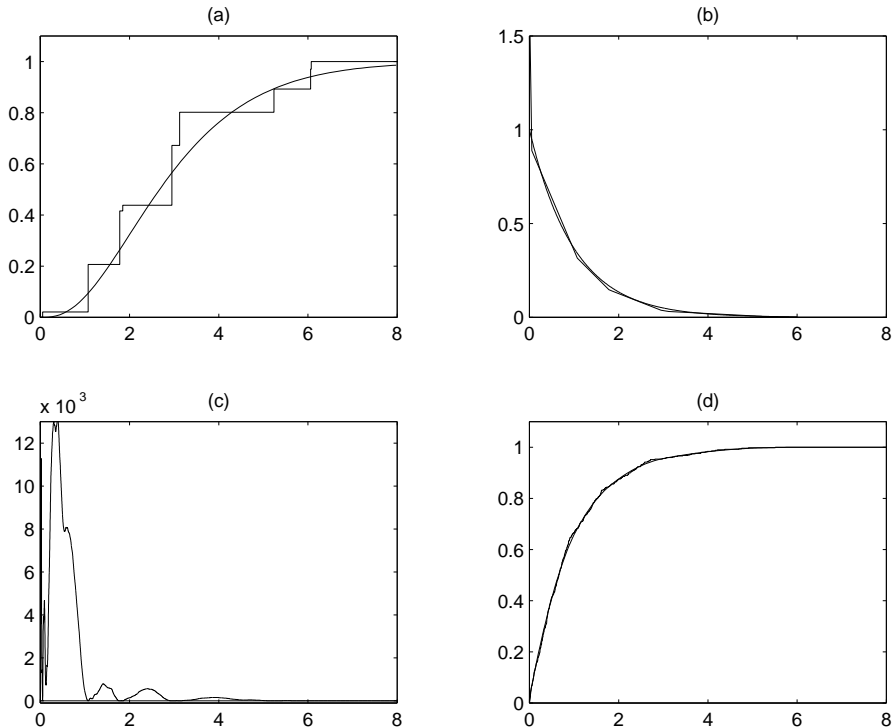}
\end{center}
\caption{{\it (a) LS estimate of the mixing distribution with
the true Gamma (3) mixing distribution; (b) LS estimate of
the mixture density with the true density; (c) the (alternative) directional
derivative function evaluated at the LS estimate and (d)
LS estimate of the mixture distribution with
the empirical distribution function of the data.
All pictures are based on a
sample of size $n=500$ from the standard exponential
distribution.}}
\label{fig:exp500CDF}
\end{figure}

\section{ML estimation in Gaussian deconvolution}
\label{sec:mldecon} In order to apply the support reduction
algorithm of section \ref{sec:simar}, the setting of Example 2 is
not appropriate since the minimization there has to be performed
over the convex hull of the functions $f_{\theta}$ instead of the
convex cone generated by them. Contrary to the situation of
section \ref{sec:lsconv}, the minimizer of $\phi$ over the cone
does not exist (given a function $f$ with $\phi(f)<0$, the
function $\phi$ applied to $c\cdot f$ for $c>0$ tends to minus
infinity). To get a well posed optimization problem over the
convex cone so that its solution is the minimizer of $\phi$ over
the convex hull of ${\cal F}$, we have to relax the constraint
that the solution has to be a probability density. The new
objective function then becomes
$$
\phi(f)=-\int \log f(x)\,d\, \FF_n(x)+\int f(x)\,dx.
$$

In principle, the support reduction algorithm can be applied
directly to the thus obtained optimization problem. However, we
observed that a Newton-type procedure (based on the support
reduction algorithm) worked significantly better than the direct
application of the support reduction algorithm. We describe this
Newton procedure here. Write $\bar f$ for the {\it current
iterate}.

Note that
$$
\phi(f)-\phi(\bar f)=-\int\log \left(1+\frac{f(x)-\bar f(x)}{\bar
f(x)}\right)\,d\FF_n(x)+ \int f(x)-\bar f(x)\,dx
$$
For $(f-\bar f)/\bar f$ small, we get the following quadratic
approximation of $\phi$ at $\bar f$, using the second order Taylor
approximation of the logarithm at $1$
$$
\int \frac12\left(\frac{f(x)-\bar f(x)}{\bar
f(x)}\right)^2-\frac{f(x)-\bar f(x)}{\bar f(x)}\,d\FF_n(x) +\int
f(x)-\bar f(x)\,dx
$$
Ignoring terms that do not depend on $f$, we define the following
local objective function
$$
\phi_q(f)=\phi_q(f;\bar f)=\int f(x)\,dx+\int
\frac12\left(\frac{f(x)}{\bar f(x)}\right)^2 -2\frac{f(x)}{\bar
f(x)}\,d\FF_n(x)
$$
This quadratic function can be minimized over the (finitely
generated) cone using the support reduction algorithm, yielding
$$
\bar{f}_q=\argmin\{\phi_q(f;\bar f)\,:\,f\in{\rm
cone}(f_{\theta}\,:\,\theta\in\Theta_{\delta})\}
$$
The next iterate is then obtained as $\bar f+\lambda(\bar f_q-\bar
f)$ ($\lambda$ chosen appropriately to assure monotonicity of the algorithm).

This method is used to solve the (finite dimensional) optimization
problem over the cone of functions generated by
$\{f_{\theta}\,:\,\theta\in\Theta_{\delta}\}$. After this, the
fine tuning in support points (leaving the prespecified grid) is
performed as described in section \ref{sec:gridless}.

During the Newton iterations to obtain the solution to the finite
dimensional problem as well as in the fine tuning step following
it,  quadratic optimization problems of the type find
$$
\argmin\{\phi_q(f)\,:\,f\in{\rm cone}(f_{\theta}\,:\,\theta\in S\}
$$
are solved for some {\it finite} set $S$. Starting from an initial
value, say $g$ (the natural candidate for this will be obvious
from the context; usually it has only a few active vertices), the
support reduction
algorithm consists of two steps that are iterated:
\begin{itemize}
    \item [1)] Find new support point
    \item [2)] Do finite dimensional constrained optimization using
iterative
unconstrained minimizations.
\end{itemize}

\vspace{0.2cm}

\noindent {\it Step 1.} In the notation of section \ref{sec:quad},
we have
$$
c_1(\theta,g)=\int
f_{\theta}(x)\,dx-2\int\frac{f_{\theta}}{\bar{f}}(x)\,d\FF_n(x)+
\int\frac{gf_{\theta}}{\bar{f}^2}(x)\,d\FF_n(x)
\mbox{ and }
c_2(\theta)=\int\frac{f_{\theta}^2}{\bar{f}^2}(x)\,d\FF_n(x).
$$
Hence, the new vertex is given by
$$
\hat{\theta}=
\argmin_{\theta\in \Theta_{\delta}} \frac{c_1(\theta,g)}{\sqrt{c_2(\theta)}}
$$

\vspace{0.2cm} \noindent {\it Step 2.} During this step, given a
support set $\{\theta_1,\ldots,\theta_p\}$, we should find a
subset $S$ of $\{\theta_1,\ldots,\theta_p\}$ with associated
optimal $f$ such that $f$ minimizes $\phi_q$ over the linear space
generated by the functions $\{f_{\theta}\,:\,\theta\in S\}$ and,
moreover, has only scalars $\alpha_j>0$ in the representation
$$
f=\sum_{\theta_j\in S}\alpha_j f_{\theta_j}
$$
The basic step in finding $S$ and $f$ is to minimize, without
restrictions on $\alpha_j$, the quadratic function
\begin{eqnarray*}
&&\psi(\alpha_1,\ldots,\alpha_p)=\phi_q(\sum_{\theta_j\in S}\alpha_j
f_{\theta_j})\\
&&\,\,\,=\sum_{i=1}^p\a_i\left(\int f_{\theta_i}(x)\,dx-2\int
\frac{f_{\theta_i}(x)}{\bar f(x)}\,d\FF_n(x)\right)
+\frac12\sum_{i=1}^p\sum_{j=1}^p\a_i\a_j\int\frac{f_{\theta_i}(x)f_{\theta_j}(x)}{\bar
f(x)^2}\,d\FF_n(x)\\
&&\,\,\,=\a^T\nu+\frac12 \a^TV\a
\end{eqnarray*}
Define the $n\times p$-matrix $Y$ by $Y_{ij}=f_{\theta_j}(x_i)$.
Note that this matrix only depends on the values of the current
vertices at the observed sample. Also define the $n$-vector $d$ by
$d_i=(\bar f(x_i))^{-1}$ and the $n\times n$ diagonal matrix $D$
$D_{ii}=d_i$. Then $nV=Y^TD^TDY$ and $n\nu=n_p-2Y^Td$ (using that
the vertices are in fact {\it probability densities}, denoting by
$n_p$ the $p$-vector with all elements equal to $n$) and the
optimal $\alpha\in\R^p$ minimizing $\psi$ is the solution to the
following linear system of equations
$$
(DY)^TDY\alpha=2Y^Td-n_p
$$
If the matrix $DY$ has full rank $p$, this system has a unique
solution.

Figure 2 shows the results of the SR algorithm based on a simulated
dataset of size $n=500$ where the mixing distribution is standard exponential.
First it took 25 iterations to obtain the solution on an equidistant grid
of size $500$ in $[x_{1)},x_{(n)}]=[-2.47,7.96]$. This grid-solution
had eight support points. After that,
1085 steps of the fine tuning step of section \ref{sec:gridless} were taken,
resulting in an estimate of the mixing distribution with five support points.

\vspace{0.2cm}

\begin{figure}[!ht]
\begin{center}
\leavevmode
\includegraphics[width=12cm]{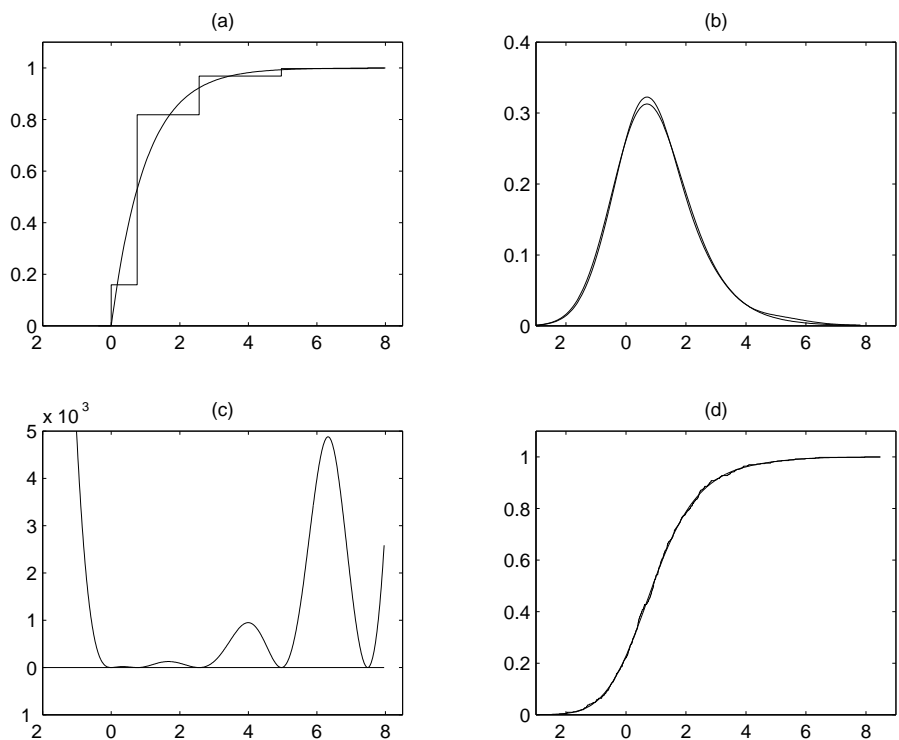}
\end{center}
\caption{{\it (a) ML estimate of the mixing distribution with
the true mixing distribution; (b) ML estimate of
the mixture density with the true density; (c) the (alternative) directional
derivative function evaluated at the ML estimate and (d)
ML estimate of the mixture distribution with
the empirical distribution function of the data.
All pictures are based on a
sample of size $n=500$ from the standard exponential
mixture of standard normals.}}
\label{fig:norm500CDF}
\end{figure}

\vspace{0.2cm}

\noindent
{\bf Acknowledgement:} We thank Jon Wakefield for drawing our attention
to Mallet's paper.






\begin{thebibliography}{}

\bibitem[B\"{o}hning (1982)]{boehning:82}
{\small \sc B\"{o}hning, D.\ (1982)}.
\newblock Convergence of Simar's algorithm for finding the maximum
likelihood estimate of a compound Poisson process.
\newblock {\sl Ann. Statist.} {\bf 10}, 1006--1008.

\bibitem[B\"{o}hning (1986)]{boehning:86}
{\small \sc B\"{o}hning, D.\ (1986)}.
\newblock A vertex-exchange method in $D$-optimal design theory.
\newblock {\sl Metrika} {\bf 33}, 337--347.

\bibitem[Dempster, Laird and Rubin (1977)]{dlr:77}
{\small \sc Dempster, A.P., Laird, N.M.\ and Rubin, D.B.\ (1977)}.
\newblock Maximum likelihood from incomplete data via the EM algorithm.
\newblock {\sl J.\ Roy.\ Statist.\ Soc.\ B} {\bf 39}, 1--38.

\bibitem[Fedorov (1972)]{fedorov}
{\small \sc Fedorov, V.V.\ (1972)}.
\newblock {\sl Theory of optimal experiments}.
\newblock Academic, New York.

\bibitem[Groeneboom, Jongbloed and Wellner (2001a)]{grojowe1:01}
{\small \sc Groeneboom, P., Jongbloed, G., and Wellner, J.A.\
(2001a)}.
\newblock A canonical process for estimation of convex functions:
the "invelope" of integrated Brownian motion $+t^4$.
\newblock {\sl  Ann.\ Statist.} {\bf 29}, 1620--1652.

\bibitem[Groeneboom, Jongbloed and Wellner (2001b)]{grojowe2:01}
{\small \sc Groeneboom, P., Jongbloed, G., and Wellner, J.A.\
(2001b)}.
\newblock Estimation of convex functions:
characterizations and asymptotic theory.
\newblock {\sl Ann.\ Statist.} {\bf 29}, 1653--1698.

\bibitem[Groeneboom and Wellner (1992)]{growe:92}
{\small \sc Groeneboom, P.\ and Wellner, J.A.\ (1992)}.
\newblock {\sl Information bounds and nonparametric
maximum likelihood estimation.}
\newblock Birkh\"auser, Basel.

\bibitem[Hampel (1987)]{hamp:87}
{\small \sc Hampel, F.R.\ (1987)}. Design, modelling and anlysis
of some biological datasets.
\newblock In {\sl Design, data and analysis, by some friends of Cuthbert
Daniel}, C.L. Mallows, editor, 111- 115.  Wiley, New York.

\bibitem[Jongbloed (1998)] {jo:98}
{\small \sc Jongbloed, G.\ (1998).}
\newblock The iterative convex minorant algorithm for nonparametric
estimation.
\newblock {\sl J.\ Comp.\ Graph.\ Statist.} {\bf 7}, 310--321.

\bibitem[Lavee, Safrie, and Meilijson (1991)]{lavsafmei:91}
{\small \sc Lavee, D., Safrie, U.N., and Meilijson, I.\ (1991)}.
\newblock For how long do trans-Sahran migrants stop over at an oasis?
\newblock {\sl Ornis Scandinavica} {\bf 22}, 33-44.

\bibitem[Lesperance and Kalbfleisch (1992)] {lk:92}
{\small \sc Lesperance, M.L.\ and Kalbfleisch, J.D.\ (1992)}.
\newblock An algorithm for computing the nonparametric MLE of a mixing
distribution.
\newblock {\sl Journal of the Americal Statistical Association} {\bf 87}
120--126.


\bibitem[Lindsay (1995)] {lindsay:95}
{\small \sc Lindsay, B.G.\ (1995)}.
\newblock {\sl Mixture models: theory, geometry and applications.}
\newblock NSF-CBMS Regional Conference Series in Probability and
Statistics, {\bf Vol 5}.

\bibitem[Mallet (1986)]{mallet:86}
{\small \sc Mallet, A.\ (1986)}.
\newblock A maximum likelihood estimation method for random coefficient
regression models.
\newblock {\sl Biometrika} {\bf 73}, 645--656.


\bibitem[Robertson, Wright, and Dykstra (1988)]{robwridyk:88}
{\small \sc Robertson, T., Wright, F. T., Dykstra, R. L.\ (1988)}.
\newblock {\sl Order Restricted Statistical Inference.}
\newblock Wiley, New York.

\bibitem[Simar (1976)]{simar:76}
{\small \sc Simar, L.\ (1976)}.
\newblock Maximum likelihood estimation of a compound Poisson process.
\newblock {\sl Ann.\ Statist.} {\bf 4}, 1200--1209.


\bibitem[Wellner and Zhang (2000)]{wezha:00}
{\small \sc Wellner, J.A.\ and Zhang, Y.\ (2000)}.
\newblock Two estimators of the mean of a counting process with panel
count data.
\newblock {\sl Ann.\ Statist.} {\bf 28}, 779--814.

\bibitem[Wu (1978)]{wu:78}
{\small \sc Wu, C.F.\ (1978)}.
\newblock Some algorithmic aspects of the theory of optimal design.
\newblock {\sl Ann.\ Statist.} {\bf 6}, 1286--1301.

\bibitem[Wynn (1970)]{wynn}
{\small \sc Wynn, H.P.\ (1970)}.
\newblock The sequential generation of $D$-optimum experimental designs.
\newblock {\sl Ann.\ Math.\ Statist.} {\bf 6} 1286--1301.

\end{thebibliography}
\end{document}